\documentclass[11pt]{article}
\usepackage{epsfig}
\usepackage{amssymb}
\usepackage{amscd}
\usepackage[all]{xy}
\usepackage{epsf}
\usepackage{amsfonts,amssymb,amsthm,amsmath}

\input cyracc.def

\long\def\comment#1\endcomment{\relax}

\newcounter{subsubsubsection}
\newcounter{subsubsubsubsection}

\makeatletter \@addtoreset{subsubsubsection}{subsubsection}
\@addtoreset{subsubsubsubsection}{subsubsubsection}

\makeatother

\DeclareMathSizes{11.1}{10}{8}{6}

\DeclareMathOperator{\Hom}{Hom}

\newtheorem{theorem}{Theorem}

\newtheorem*{theorem*}{Theorem}

\newtheorem*{lemma}{Lemma}

\newtheorem*{proposition}{Proposition}

\newtheorem*{corollary}{Corollary}

\theoremstyle{remark}

\newtheorem*{remark}{Remark}

\theoremstyle{definition}

\DeclareMathOperator{\Id}{Id}

\DeclareMathOperator{\ad}{ad}

\newcommand{\hdot}{{\:\protect\raisebox{3pt}{\text{\protect\circle*{1.5}}}}}

\newcommand{\mb}{\hdot}

\newcommand\Cyl{\mathrm{Cyl}}

\newcommand{\g}{\mathfrak{g}}

\newcommand{\F}{\mathcal{F}}

\title{ {\huge  An explicit construction of the Quillen homotopical category of dg Lie algebras}}
\author{{\LARGE Boris Shoikhet}}
\date{}

\begin{document}\maketitle

\comment

{\font\tcyr=wncyi10

\tcyr\cyracc

\font\scyr=wncyr10

\scyr\cyracc

\hbox to\textwidth{\hfil\parbox{80mm}{\tcyr{Zachem ya algebry ne
znal podi po{\u\i}mi se{\u\i}chas}\\{\tcyr Tako{\u\i} polezny{\u\i}
delovo{\u\i} predmet...}}} \vspace{2mm}

\hbox to\textwidth{\hfil\parbox{30mm}{\scyr{M.Shcherbakov}}}}

\vspace{1cm}

\endcomment
\begin{abstract}
Let $\g_1$ and $\g_2$ be two dg Lie algebras, then it is well-known
that the  $L_\infty$ morphisms from $\g_1$ to $\g_2$ are in $1-1$
correspondence to the solutions of the Maurer-Cartan equation in
some dg Lie algebra $\Bbbk(\g_1,\g_2)$. Then the gauge action by
exponents of the zero degree component $\Bbbk(\g_1,\g_2)^0$ on
$MC\subset\Bbbk(\g_1,\g_2)^1$ gives an explicit "homotopy relation"
between two $L_\infty$ morphisms. We prove that the quotient
category by this relation (that is, the category whose objects are
$L_\infty$ algebras and morphisms are $L_\infty$ morphisms modulo
the gauge relation) is well-defined, and is a localization of the
category of dg Lie algebras and dg Lie maps by quasi-isomorphisms.
As localization is unique up to an equivalence, it is equivalent to
the Quillen-Hinich homotopical category of dg Lie algebras [Q1,2],
[H1,2]. Moreover, we prove that the Quillen's concept of a homotopy
coincides with ours. The last result was conjectured by V.Dolgushev
[D].
\end{abstract}

\section{Introduction}
\subsection{}
Let $\mathcal{C}$ be a small category, and let $S\subset
\mathrm{Mor}(\mathcal{C})$ be a set of morphisms. In many problems
it is useful to construct a localized category $S^{-1}\mathcal{C}$
and a functor $P_S\colon\mathcal{C}\to S^{-1}\mathcal{C}$ which obey
the following properties:
\begin{itemize}
\item[1.] for any morphism $s\in S$, the morphism $P_S(s)$ is
invertible,
\item[2.] any other functor $F\colon \mathcal{C}\to\mathfrak{X}$
which maps any morphism $s\in S$ to an invertible morphism, can be
decomposed $F=G\circ P_S$ where the functor $G\colon
S^{-1}\mathcal{C}\to\mathfrak{X}$ is uniquely defined.
\end{itemize}
If the category $S^{-1}\mathcal{C}$ exists, it is unique up to an
equivalence of categories.

For any $S$, there is a general construction of $S^{-1}\mathcal{C}$,
which is a big category. This construction goes as follows.

Recall that a diagram scheme $T$ is a set of objects
$\mathrm{Ob}(T)$ and a set of morphisms $\mathrm{Mor}(T)$ with no
compositions. For each morphism $m\in \mathrm{Mor}(T)$ there defined
its beginning $\alpha(m)$ and its end $\beta(m)$. One can attach a
diagram scheme to each category forgetting the compositions.

For each diagram scheme, one associates a category of paths
$\mathcal{P}(T)$ as follows:
$\mathrm{Ob}(\mathcal{P}(T))=\mathrm{Ob}(T)$, and the morphisms
$\mathrm{Mor}(A,B)$ are all paths, that is all the sequences
$t_1,\dots, t_k\in\mathrm{Mor}(T)$ such that $\alpha(t_1)=A$,
$\beta(t_k)=B$, and for $1\le i<k$ one has
$\alpha(t_{i+1})=\beta(t_i)$. The composition is the composition of
paths.

With a small category $\mathcal{C}$ and a set
$S\subset\mathrm{Mor}\mathcal{C}$, one associates a diagram scheme
$T$ with the set of objects $\mathrm{Ob}\mathcal{C}$, and the set
$\mathrm{Mor}(T)=\mathrm{Mor}(\mathcal{C})\sqcup S$. We denote by
$i\colon \mathrm{Mor}(\mathcal{C})\to \mathrm{Mor}(T)$ and $j\colon
S\to \mathrm{Mor}(T)$ the natural imbeddings to the first and to the
second components. For the first component, we define $\alpha(m)$
and $\beta(m)$ as the beginning and the end in the category
$\mathcal{C}$, and for the second component we define $\alpha(s)$ as
the {\it end} of $s$ in $\mathcal{C}$, and $\beta(s)$ as the {\it
beginning} of $s$ in $\mathcal{C}$. Then we have the category of
paths $\mathcal{P}(T)$.

By definition, the category $S^{-1}\mathcal{C}$ has the same objects
as $\mathcal{P}(T)$ (and, therefore, the same as $\mathcal{C})$, and
the morphisms is the quotient of the morphisms in $\mathcal{P}(T)$
by the following relations:
\begin{itemize}
\item[(i)] $i(m_1)\circ i(m_2)=i(m_1\circ m_2)$ if the composition
$m_1\circ m_2$ in $\mathcal{C}$ is defined,
\item[(ii)] $\Id_{\mathcal{C}}(a)=\Id_{\mathcal{P}(T)}(a)$ for any
object $a$,
\item[(iii)] $i(s)\circ j(s)=\Id$, $j(s)\circ i(s)=\Id$, for any
$s\in S$.
\end{itemize}

One easily sees that the conditions 1 and 2 above are satisfied.

The lack of this construction is that in this way in general we get
a category with the morphisms forming as a set a higher universe
than the morphisms in $\mathcal{C}$, that is, a big category. In
particular, it is impossible to answer any direct question, for
example, is a diagram commutative or not.

Sometimes the category  $S^{-1}\mathcal{C}$ is small. It can be seen
only by an alternative more direct construction. The most known case
is the construction of the derived category of an Abelian category.
In this example the category $\mathcal{C}$ is the category of
complexes in the Abelian category, and the set $S$ is the set of
{\it quasi-isomorphisms}, that is, maps of complexes inducing
isomorphism on cohomology. This construction uses essentially that
the set of morphisms between any two objects is an Abelian group (or
a vector space).

\subsection{}
There are many examples in which the set of morphisms is {\it not}
an Abelian group, like in the category of topological spaces, where
$S$ could be the set of homotopical equivalences.  As more advanced
example, one can consider the category of associative (or
commutative) dg algebras, where $S$ is the set of maps of algebras
which are quasi-isomorphisms of complexes. More generally, one can
consider the category of dg modules over an operad, with $S$ equal
to the set of quasi-isomorphisms.

The problem is the same--to construct explicitly the category
$S^{-1}\mathcal{C}$ as a small category, such that a question of the
commutativity of diagrams can be effectively solved. This is done by
Quillen [Q1,2], who introduced concept of a closed model category.
We recall his construction in Section 3. The idea is to axiomatize
some data in the category $\mathcal{C}$ and in $S$, such that the
quotient category could be considered as the homotopical category of
topological space. The problem is to derive the concept of homotopy,
such that the quotient category with the same objects and the
morphisms equal to the quotient sets by the homotopy relation is a
localization, in particular, the morphisms in $S$ are invertible in
the quotient category. Note that this homotopy relation it is not
easy to derive, it does not follow directly from the set $S$, and
the possibility to deal with dg algebras as with topological spaces
was a remarkable invention.

In the Quillen construction, this homotopy relation is not very
explicit. Our goal in this paper is to derive this relation in the
category of dg Lie algebras explicitly. Our construction can be
generalized for the category of dg modules over any Koszul operad.
Note that we should work with unbounded complexes, while Quillen
[Q1,2] considered only bounded complexes. For the unbounded case,
the construction was generalized by Hinich [H1,2].

\subsection{}
The paper is organized as follows.

Section 2 contains a definition of the homotopy relation between two
$L_\infty$ morphisms for dg Lie algebras; here we prove that the
quotient category $\mathcal{H}om_{dg}$ by this homotopy relation is
well-defined,

Section 3 contains some straightforward generalizations of some
results of Section 2 to the case of $L_\infty$ algebras,

In Section 4 we prove that if an $L_\infty$ morphism between two dg
Lie algebras is an $L_\infty$ quasi-isomorphism, it is homotopically
invertible, that is, invertible in the category
$\mathcal{H}om_{dg}$; we also note here that for general $L_\infty$
algebras this is (probably) not true,

Section 5 contains a proof of the universal property (2) of Section
1.1 for the category of dg Lie algebras and quasi-isomorphisms; the
results of Sections 4 and 5 together show that the category
$\mathcal{H}om_{dg}$ is a localization of the category of dg Lie
algebras and dg Lie maps by quasi-isomorphisms. As the localization
is unique up to an equivalence, the category $\mathcal{H}om_{dg}$
coincides with the Quillen homotopical category; moreover, we prove
that our relation of homotopy coincides with the Quillen's one.

The paper is completely independent on closed model categories, any
knowledge on them is not supposed. We tried to write a paper
understandable for a reader without any special background.
\subsection{}
The author is grateful to Borya Feigin, Volodya Hinich,  Dima
Kaledin, Bernhard Keller, and to Dima Tamarkin for many very
interesting and useful discussions. The work was partially supported
by the research grant R1F105L15 of the University of Luxembourg.

\section{The homotopy relation}
Consider the category $\mathcal{L}_{dg}$ of dg Lie algebras and
$L_\infty$ maps. In this Section we introduce a homotopy relation
between two morphisms in this category.

Let $\g_1$ and $\g_2$ be two dg Lie algebras. Here we recall (see,
e.g. [H1], [Kel2], [D]) the construction of a dg Lie algebra
$\Bbbk(\g_1,\g_2)$ such that the solutions of the Maurer-Cartan
equation in $\Bbbk^1(\g_1,\g_2)$ are in $1-1$ correspondence with
the $L_\infty$ morphisms from $\g_1$ to $\g_2$.

As a dg vector space, $$
\Bbbk(\g_1,\g_2)=\Hom(C_+(\g_1,\mathbb{C}),\g_2)
$$
where homomorphisms are taken over the ground field (we always
suppose that it is the field of complex numbers). Here
$C(\g_1,\mathbb{C})$ is the chain complex of the dg Lie algebra
$\g_1$, it is naturally a counital dg coalgebra, and
$C_+(\g_1,\mathbb{C})$ is the kernel of the counit map. When $\g_1$
is usual Lie algebra concentrated in degree 0,
$C_+(\g_1,\mathbb{C})$ is $\mathbb{Z}_{<0}$-graded dg coalgebra.
Define now a Lie bracket on $\Bbbk(\g_1,\g_2)$. Let
$\theta_1,\theta_2\in\Bbbk(\g_1,\g_2)$ be two elements. Their
bracket $[\theta_1,\theta_2]$ is defined (up to a sign specified
below) as
\begin{equation}\label{eq2.1}
C_+(\g_1,\mathbb{C})\xrightarrow{\Delta}C_+(\g_1,\mathbb{C})^{\otimes
2}\xrightarrow{\theta_1\otimes\theta_2}\g_2\otimes\g_2\xrightarrow{[,]}g_2
\end{equation}
where $\Delta$ is the coproduct in $C_+(\g_1,\mathbb{C})$ and $[,]$
is the Lie bracket in $\g_2$. It follows from the cocommutativity of
$\Delta$ that in this way we get a Lie algebra. When $\g_1$ is
finite-dimensional, the bracket reduces to the usual bracket
\begin{equation}\label{signs}
[m_1\otimes g_1,m_2\otimes g_2]=(-1)^{\deg g_1\cdot\deg
m_2}(m_1m_2)\otimes [g_1,\g_2]
\end{equation}
on the product $m\otimes \g$ of a (graded) commutative dg algebra
$m$ with a Lie algebra $\g$. The sign rule (\ref{signs}) guarantees
that if one uses the Koszul sign rules in $\g$, that is,
$[g_1,g_2]=(-1)^{\deg g_1\deg g_2+1}[g_2,g_1]$, then $m\otimes\g$
also obeys the Koszul sign rule.

An element $F$ of degree 1 in $\Bbbk(\g_1,\g_2)$ is a collection of
maps
\begin{equation}\label{eq2.2}
\begin{aligned}
\ & F_1\colon \g_1\to\g_2\\
&F_2\colon \Lambda^2(\g_1)\to\g_2[-1]\\
&F_3\colon \Lambda^3(\g_1)\to\g_2[-2]\\
&\dots
\end{aligned}
\end{equation}
and the Maurer-Cartan equation $d_{\Bbbk}F+\frac12[F,F]_{\Bbbk}=0$
is the same that the collection $\{F_i\}$ are the Taylor components
of an $L_\infty$ map which we denote also by $F$. Note that the
differential in $\Bbbk(\g_1,\g_2)$ comes from 3 differentials: the
both inner differentials in $\g_1$ and $\g_2$, and from the chain
differential in $C_+(\g_1,\mathbb{C})$.

Now for any dg Lie algebra $\g$, the solutions of the Maurer-Cartan
equation form a quadric in $\g^1$, and $\g^0$ acts on this quadric
by vector fields. Namely, each $X\in\g^0$ defines a vector field
\begin{equation}\label{eq2.3}
\frac{dF}{dt}=-dX+[X,F]
\end{equation}
It can be directly checked that this vector field indeed preserves
the quadric.

In our case, this vector field can be exponentiated to an action on
the pro-nilpotent completion on $\Bbbk$. This action gives our
homotopy relation on $L_\infty$ morphisms.

More precisely, an element $H$ of degree 0 in $\Bbbk$ is a
collection of maps

\begin{equation}\label{eq2.4}
\begin{aligned}
\ &H_1\colon \g_1\to \g_2[-1]\\
&H_2\colon\Lambda^2(\g_1)\to\g_2[-2]\\
&H_3\colon\Lambda^3(\g_1)\to\g_2[-3]\\
&\dots
\end{aligned}
\end{equation}
We can rewrite (\ref{eq2.3}) as:
\begin{equation}\label{eq2.5}
\frac{d(F+d_\Bbbk)}{dt}=[H,F+d_\Bbbk]
\end{equation}
Then
\begin{equation}\label{eq2.6}
F_H+d_\Bbbk=\exp(ad(H))(F+d_\Bbbk)
\end{equation}
from where we find an explicit formula:
\begin{equation}\label{eq2.6}
F_H=\exp(ad(H))(F)+f(H)(d_\Bbbk H)
\end{equation}
where $f(z)=(\exp(z)-1)/z$.

We say that two $L_\infty$ morphisms $F_1,F_2$ from $\g_1$ to $\g_2$
are {\it homotopic} if $F_2=(F_1)_H$ for some $H$ as above.

It is clear that it is an equivalence relation. Namely, the
reflexivity is clear, as well as the symmetry. To prove the
transitivity, we should use the Campbell-Baker-Hausdorff formula,
and it goes without problems.

We are going to prove that this relation is compatible with the {\it
composition} of $L_\infty$ morphisms.

Firstly recall the following lemma:
\begin{lemma}
Let $\g_1$ and $\g_2$ be two differential graded Lie algebras, and
$\mathcal{U}\colon\g_1\to\g_2$ be an $L_\infty$ morphism. Then for
any solution $\gamma\in\g_1^1$ of the Maurer-Cartan equation the
$L_\infty$ morphism $\mathcal{U}$ defines a solution
$\mathcal{U}_*(\gamma)\in\g_2^1$ of the Maurer-Cartan equation in
$\g_2$. Moreover, if two solutions $\gamma_1$ and $\gamma_2$
obtained one from another by the exponentiated action of
$x\in\g_1^0$, the solutions $\mathcal{U}_*(\gamma_1)$ and
$\mathcal{U}_*(\gamma_2)$ are also obtained one from another by the
exponentiated action of an element $\mathcal{U}_*(x)\in\g_2^0$. That
is, the map $\mathcal{U}_*$ maps gauge equivalent solutions of the
Maurer-Cartan equation to gauge equivalent ones.
\begin{proof}
An $L_\infty$ algebra structure on $\g$ is by definition the same
that a vector field $Q$ of degree +1 on the space $\g[1]$ such that
$[Q,Q]=0$. One can speak on zeros of this odd field $Q$. Namely, in
usual situation a vector field $v$ on some space $X$ vanishes in a
point $p\in X$ iff for any function $f$ in the neighborhood of $p$
the value $v(f)$ vanishes at $p$. The same definition works in the
odd case as well. Let $\g$ be a dg Lie algebra. One can prove that
the Maurer-Cartan quadric is the same that the zero locus of the
corresponding field $Q$. Moreover, by definition, an $L_\infty$
morphism is a $Q$-equivariant map, and therefore it maps a point on
the zero locus of $Q_{\g_1}$ to a point on the zero locus of
$Q_{\g_2}$. Thus, a solution of the Maurer-Cartan equation is mapped
by an $L_\infty$ map to a solution of the Maurer-Cartan equation.

In the Taylor components, the formula for $\mathcal{U}_*(\gamma)$ is
\begin{equation}\label{eqlemma2.1.1}
\mathcal{U}_*(\gamma)=\mathcal{U}_1(\gamma)+\frac12\mathcal{U}_2(\gamma,\gamma)+...+\frac1{n!}\mathcal{U}_n(\gamma,\dots,\gamma)+\dots
\end{equation}
By the same reasons, an $L_\infty$ map maps a vector field tangent
to the Maurer-Cartan quadric to a tangent vector field. If the first
tangent vector field is corresponded to an element $x\in\g_1^0$, the
second one is corresponded to the element
\begin{equation}\label{eqlemma2.1.2}
\mathcal{U}_*(x)=\mathcal{U}_1(x)+\mathcal{U}_2(x,\gamma)+\frac12\mathcal{U}_3(x,\gamma,\gamma)+...+\frac1{n!}\mathcal{U}_{n+1}(x,\gamma,\dots,\gamma)+\dots
\end{equation}
\end{proof}
\end{lemma}

Now we are going to prove the following proposition:

\begin{proposition}
Let $\g_1,\g_2,\g_3$ be two dg Lie algebras, and let $F\colon
\g_1\to\g_2$ and $G\colon\g_2\to\g_3$ be two $L_\infty$ morphisms.
Suppose that $F_1=F_H$ and $G_1=G_U$ for some homotopies
$H\in\Bbbk^0(\g_1,\g_2)$ and $U\in\Bbbk^0(\g_2,\g_3)$. Then
$G_1\circ F_1=(G\circ F)_X$ for some homotopy
$X\in\Bbbk^0(\g_1,\g_3)$.
\begin{proof}
The maps $F$ and $G$ define the following diagram of $L_\infty$
morphisms:
\begin{equation}\label{eq2.7}
\xymatrix{& \Bbbk(\g_1,\g_3)\\
\Bbbk(\g_2,\g_3)\ar[ur]^{F^!}&& \Bbbk(\g_1,\g_2)\ar[ul]_{G_!}}
\end{equation}
Now we consider $G$ as a solution of the Maurer-Cartan equation in
$\Bbbk(\g_2,\g_3)$ and $F$ as a solution of the Maurer-Cartan
equation in $\Bbbk(\g_1,\g_2)$. Then
\begin{equation}\label{eq2.8} (F^!)_*(G)=(G_!)_*(F)=G\circ F
\end{equation}
where $G\circ F$ is the composition of $F$ and $G$ considered as a
solution of the Maurer-Cartan equation in $\Bbbk(\g_1,\g_3)$.

Then for any homotopies $H$ and $U$ as above, it follows from the
Lemma that
$$
(F^!)_*(G_U)\sim (F^!)_*(G)=(G_!)_*(F)\sim(G_!)_*(F_H)
$$
where $\sim$ stands for the gauge equivalence of the solutions of
the Maurer-Cartan equation.

This means that $G_U\circ F\sim G\circ F_H$ which immediately
implies the statement of Proposition.
\end{proof}
\end{proposition}

The Proposition implies that the homotopical equivalence between
$L_\infty$ morphisms is compatible with the composition. We define
the category $\mathcal{H}om_{dg}$ as follows: its objects are
differential graded Lie algebras, and its morphisms are $L_\infty$
morphisms modulo the homotopy relation.

We prove in Section 2.3 the following theorem:

\begin{theorem}
Any $L_\infty$ quasi-isomorphism between two dg Lie algebras is
invertible in the category $\mathcal{H}om_{dg}$
\end{theorem}

It is clear that if an $L_\infty$ morphism is invertible in
$\mathcal{H}om_{dg}$, it is an $L_\infty$ quasiisomorphism.

Before proving the theorem, we generalize some of the constructions
of this Section to the case of $L_\infty$ algebras.

\section{The case of $L_\infty$ algebras}
Recall that $\g$ is an $L_\infty$ algebra if the cofree
cocommutative dg coalgebra without counit $Fun_+(\g[1])$ is endowed
with a coderivation $Q$ of degree $+1$ such that $[Q,Q]=0$. In this
case we still denote the complex $Fun_+(\g[1],Q)$ by
$C_+(\g,\mathbb{C})$.

If $\g_1,\g_2$ are $L_\infty$ algebras, then
$\Bbbk(\g_1,\g_2)=\Hom(C_+(\g_1),\g_2)$ is again an $L_\infty$
algebra. The Taylor components of this structure
$$
\mathcal{L}_k\colon\Lambda^k(\Bbbk(\g_1,\g_2))\to\Bbbk(\g_1,\g_2)[2-k]
$$
are defined as follows:

The component $\mathcal{L}_0$ (the differential) comes from the
$L_\infty$ differential $Q$ on $C_+(\g_1,\mathbb{C})$ and from the
inner differential on dg space $\g_2$ (that is, from the first
Taylor component of the $L_\infty$ structure on $\g_2$).

For the higher component we have:

To define $\mathcal{L}_k$, we first apply the $(k-1)$-st power of
the coproduct in $C_+(\g_1,\mathbb{C})$ and obtain an element in
$C_+(\g_1,\mathbb{C})^{\otimes k}$, then apply
$\Psi_1\otimes\dots\otimes \Psi_k$ where
$\Psi_i\in\Bbbk(\g_1,\g_2)$, we get an element in $\otimes^k(\g_2)$,
and then apply the $k$-th Taylor component of the $L_\infty$
structure on $\g_2$. More precisely, it is
$$
C_+(\g_1,\mathbb{C})\xrightarrow{\Delta^{k-1}}C_+(\g_1,\mathbb{C})^{\otimes
k}\xrightarrow{\Psi_1\otimes\dots\otimes\Psi_k}\g_2^{\otimes
k}\xrightarrow{\mathcal{L}_k(\g_2)}\g_2
$$
where $\mathcal{L}_k(\g_2)$ is the $k$-th Taylor component of the
$L_\infty$ structure on $\g_2$.

For an $L_\infty$ algebra $\g$, one can define a solution of the
Maurer-Cartan equation as $\alpha\in\g^1$ such that
\begin{equation}\label{eqprooff1}
\mathcal{L}_1(\alpha)+\frac12\mathcal{L}_2(\alpha,\alpha)+\dots+\frac1{k!}\mathcal{L}_k(\alpha,\dots,\alpha)+\dots=0
\end{equation}
where $\mathcal{L}_k$'s are the Taylor components of the $L_\infty$
structure on $\g$.

Also, an element $x\in\g^0$ acts on the solutions of the
Maurer-Cartan equation by the formula
\begin{equation}\label{eqprooff2}
\frac{d\alpha}{dt}=-\mathcal{L}_1(x)+\mathcal{L}_2(x,\alpha)-\frac12\mathcal{L}_3(x,\alpha,\alpha)+\dots
\end{equation}

One can prove that the solutions of the Maurer-Cartan equation in
$\Bbbk(\g_1,\g_2)$ where $\g_1,\g_2$ are two $L_\infty$ algebras,
are exactly the $L_\infty$ morphisms from $\g_1$ to $\g_2$.

Then the space $\Bbbk(\g_1,\g_2)^0$ acts on the solutions of the
Maurer-Cartan equation, as in the case of dg Lie algebras in Section
2.

In the case $\g=\Bbbk(\g_1,\g_2)$ the exponent of such a vector
field is well-defined, and gives us the definition of homotopic
$L_\infty$ morphisms.

If two $L_\infty$ algebras are $L_\infty$ quasi-isomorphic, then the
solutions of the Maurer-Cartan equation modulo this gauge action
define equivalent deformation functors.

The Lemma and Proposition in Section 2 are true in the $L_\infty$
case, and we can define the category $\mathcal{H}om_\infty$ whose
objects are $L_\infty$ algebras and the morphisms between $L_\infty$
algebras $\g_1$ and $\g_2$ are the $L_\infty$ morphisms modulo the
gauge equivalence, obtained by the integration of the vector fields
(\ref{eqprooff2}) on the Maurer-Cartan quadric in
$\Bbbk(\g_1,\g_2)^1$.

Let us note that it is {\it not} true in general that an $L_\infty$
quasi-isomorphism between two $L_\infty$ algebras has a
homotopically inverse. That is, the analog of Theorem 1 for
$L_\infty$ algebras fails. See some explanation of this phenomenum
in the remark after Lemma 4.3.

\section{A proof of Theorem 1}
The proof of Theorem 1 is divided by several steps.
\subsection{}
Let $\g_1,\g_2$ be two dg algebras, and let $\mathcal{F}\colon
\g_1\to\g_2$ be an $L_\infty$ quasi-isomorphism. Using the
construction from [Me],[P],[KS] with planar trees (which is a more
direct version of the Kadeishvili theorem), we have the "induced"
$L_\infty$ structure on the cohomology (of complexes) $H^\mb(\g_1)$
and $H^\mb(\g_2)$ which are $L_\infty$ quasi-isomorphic to $\g_1$
and $\g_2$, correspondingly. There are $L_\infty$ quasi-isomorphisms
$\mathcal{M}_i\colon H^\mb(\g_i)\to\g_i$, $(i=1,2)$, which can be
constructed also by planar trees (see [KS]). We have the following
solid arrow diagram:
\begin{equation}\label{eq4.1}
\xymatrix{ \g_1\ar[rr]^{\mathcal{F}}&&\g_2\\
\\
H^\mb(\g_1)\ar[uu]^{\mathcal{M}_1}\ar@{-->}[rr]&&H^\mb(\g_2)\ar[uu]_{\mathcal{M}_2}}
\end{equation}
The $L_\infty$ structures on $H^\mb(\g_i)$ and the maps
$\mathcal{M}_i$ are not constructed canonically, they depend on a
splitting of the {\it complexes} $\g_i=H^\mb(\g_i)\oplus L_i$ where
$L_i$ are acyclic subcomplexes. In particular, the dotted arrow at
the moment is not constructed.

We need to know only one thing about this planar trees construction
of "Massey operations", which we fix in a separate lemma:

\begin{lemma}
The maps $\mathcal{M}_i\colon H^\mb(\g_i)\to\g_i$ are $L_\infty$
imbeddings, that is, the first Taylor component of the $L_\infty$
morphisms $\mathcal{M}_i$ are imbeddings of complexes. Moreover,
$(\mathcal{M}_i)_1$ is a map of complexes which induces the identity
isomorphism on the cohomology.
\begin{proof}
It follows from the construction with the splitting
$\g_i=H^\mb(\g_i)\oplus L_i$ and with planar trees, see [KS] for
details.
\end{proof}
\end{lemma}
\subsection{}
\begin{lemma}
Let $\mathcal{M}\colon\mathfrak{t}_1\to\mathfrak{t}_2$ be an
$L_\infty$ imbedding (see the definition in Lemma 4.1). Then the
corresponding map of cofree cocommutative dg coalgebras
$\mathcal{M}_*\colon C_+(\mathfrak{t}_1,\mathbb{C})\to
C_+(\mathfrak{t}_2,\mathbb{C})$ is also injective map of vector
spaces.
\begin{proof}
Let $T\in C_+(\mathfrak{t}_1,\mathbb{C})$ belongs to the kernel of
$\mathcal{M}_*$. Then $T$ is a finite sum, that is, $T\in
\oplus_{i\le k}\Lambda^i(\mathfrak{t}_1)$. Consider the highest
degree part $T_k\in \Lambda^k(\mathfrak{t}_1)$. Then
\begin{equation}\label{eqproof2}
\mathcal{M}_*(\Delta^{k-1}(T_k))=\mathcal{M}(\Delta^{k-1}(T))=\Delta^{k-1}(\mathcal{M}(T))=0
\end{equation}
But $\Delta^{k-1}(T_k)\in \Lambda^k(\mathfrak{t}_1)$. Then
$\mathcal{M}$ acts component-wise and on $\mathfrak{t}_1$ it acts
injectively. This proves that $T_k=0$, and the assertion of lemma.
\end{proof}
\end{lemma}
Clearly this lemma has a more intuitively evident counterpart for
algebras, when the injectivity is replaced by the surjectivity.

\subsection{The Quillen functors}
Up to now, we reduced Theorem 1 to the case of $L_\infty$
quasi-isomorphisms which are also imbeddings, now we make one more
reduction and reduce the Theorem to the case of quasi-isomorphic
imbeddings of dg Lie algebras.

Let $\mathsf{Lie}$ be the category of (unbounded) dg Lie algebras
over $\mathbb{C}$ with maps of dg Lie algebras as morphisms, and let
$\mathsf{Coalg}$ be the category of (unbounded) counital
cocommutative dg coalgebras over $\mathbb{C}$ with maps of
coalgebras as morphisms. There are two functors
$$
\mathcal{L}\colon \mathsf{Coalg}\to\mathsf{Lie}
$$
and
$$
\mathcal{C}\colon \mathsf{Lie}\to\mathsf{Coalg}
$$
such that $\mathcal{L}$ is the left adjoint to $\mathcal{C}$, and
$\mathcal{C}$ is the right adjoint to $\mathcal{L}$. That is,
\begin{equation}\label{eqproof25}
\Hom_{\mathsf{Lie}}(\mathcal{L}(X),Y)=\Hom_{\mathsf{Coalg}}(X,\mathcal{C}(Y))
\end{equation}
Let $X$ be a counital cocommutative dg coalgebra, and let
$\bar{X}=\mathrm{Ker}\{\varepsilon\colon X\to\mathbb{C}\}$ be the
kernel of the counit. Consider the free Lie algebra generated by
$\bar{X}[-1]$ endowed with a differential arose from the coproduct
$\Delta\colon
\bar{X}[-1]\to\Lambda^2(\bar{X}[-1])[1]=S^2(\bar{X})[-1]$. This is
the dg Lie algebra $\mathcal{L}(X)$. For the functor $\mathcal{C}$,
if $Y$ is a dg Lie algebra, $\mathcal{C}(Y)$ is the chain complex
$C(Y,\mathbb{C})=S(Y[1])$, with the natural dg coalgebra structure.

From the adjointness property (\ref{eqproof25}), we have natural
adjunction morphisms
$$
\varepsilon_{\mathsf{Coalg}}\colon X\to \mathcal{C}\circ\mathcal{L}
(X)
$$
and
$$
\varepsilon_{\mathsf{Lie}}\colon \mathcal{L}\circ\mathcal{C}(Y)\to Y
$$
which are both quasi-isomorphisms.

The Quillen functors allow to reduce a question about an $L_\infty$
map from an $L_\infty$ algebra to dg Lie algebra to a dg Lie map
between two dg Lie algebras. We have the following lemma:

\begin{lemma}
Suppose that $\mathfrak{a}_1$ is an $L_\infty$ algebra, and
$\mathfrak{a}_2$ is a dg Lie algebra. Then the correspondence
$Q\colon
\mathcal{M}\rightsquigarrow\mathcal{G}=\mathcal{L}(\mathcal{M})$
from $L_\infty$ morphisms from $\mathfrak{a}_1$ to $\mathfrak{a}_2$
to dg Lie algebras maps from
$\mathcal{L}(C_+(\mathfrak{a}_1,\mathbb{C}))$ to
$\mathcal{L}(C_+(\mathfrak{a}_2,\mathbb{C}))$ is a functorial $1-1$
correspondence preserving the classes of quasi-isomorphisms and of
imbeddings. If $Q(\mathcal{M}_1)$ and $Q(\mathcal{M}_2)$ are
homotopic, then $\mathcal{M}_1$ and $\mathcal{M}_2$ also are
homotopic.
\begin{proof}
Let $\mathcal{G}\colon\mathcal{L}(C_+(\mathfrak{a}_1,\mathbb{C}))\to
\mathcal{L}(C_+(\mathfrak{a}_2,\mathbb{C}))$ be a map of dg Lie
algebras. We want to reconstruct the map $\mathcal{M}$.

Consider the composition
\begin{equation}\label{eqnew4.1}
C_+(\mathfrak{a}_1,\mathbb{C})\xrightarrow{\varepsilon_{\mathsf{Coalg}}}
\mathcal{C}\circ\mathcal{L}(C_+(\mathfrak{a}_1,\mathbb{C}))\xrightarrow
{\mathcal{C}(\mathcal{G})}\mathcal{C}\circ\mathcal{L}(C_+(\mathfrak{a}_2,\mathbb{C}))\xrightarrow
{\mathcal{C}\circ\mathcal{L}(i)}\mathcal{C}\circ\mathcal{L}\circ\mathcal{C}(\mathfrak{a}_2)\xrightarrow
{\mathcal{C}\circ\varepsilon_{\mathsf{Lie}}}\mathcal{C}(\mathfrak{a}_2)
\end{equation}
Here $i\colon C_+(\mathfrak{a}_2,\mathbb{C})\to
C(\mathfrak{a}_2,\mathbb{C})=\mathcal{C}(\mathfrak{a}_2)$ is the
canonical imbedding. The composition (\ref{eqnew4.1}) takes images
actually in $C_+(\mathfrak{a}_2,\mathbb{C})$, and defines an
$L_\infty$ morphism
$\mathcal{M}\colon\mathfrak{a}_1\to\mathfrak{a}_2$. Denote the
correspondence $\mathcal{G}\rightsquigarrow\mathcal{M}$ by $Q^{-1}$.
It is easy to prove that $Q$ and $Q^{-1}$ are inverse to each other.
Their functoriality, and preserving of the quasi-isomorphisms and of
the imbeddings, is clear. It remains to prove the statement about
homotopy.

For this we construct a map of dg Lie algebras
$$
Q_*\colon\Bbbk(\mathcal{L}(C_+(\mathfrak{a}_1,\mathbb{C})),\mathcal{L}(C_+(\mathfrak{a}_2,\mathbb{C})))\to\Bbbk(\mathfrak{a}_1,\mathfrak{a}_2)
$$
exactly in the same way as we constructed the correspondence
$Q^{-1}$ in (\ref{eqnew4.1}). It is a map of dg Lie algebras,
therefore, it maps gauge equivalent solutions of the Maurer-Cartan
equation to gauge equivalent ones.
\end{proof}
\end{lemma}

\begin{remark}
The statement of lemma is not true when $\mathfrak{a}_2$ is an
$L_\infty$ algebra but not a dg Lie algebra. It shows that an
$L_\infty$ quasi-isomorphism between two $L_\infty$ algebras in
general is not homotopically invertible, when the homotopy is
understood in the sense of Section 2. It could be interesting to
find an example of such an $L_\infty$ quasi-isomorphism between two
$L_\infty$ algebras.
\end{remark}

\subsection{The case of a quasi-isomorphic imbedding of dg Lie
algebras}
\begin{proposition}
Let $i\colon \g_0\to \g$ be an imbedding of dg Lie algebras which is
a quasi-isomorphism. Then there exists an $L_\infty$ morphism
$\mathcal{F}\colon \g\to \g_0$ such that the both compositions
$\mathcal{F}\circ i$ and $i\circ\mathcal{F}$ are homotopic to
identity maps.
\begin{proof}
Firstly we construct an $L_\infty$ map $\mathcal{F}\colon
\g\to\g_0$.

We split $\g=\g_0\oplus L$ where $L$ is an acyclic {\it complex}.
Such a splitting always exists, but is not canonical. Then we have a
projection $p\colon \g\to\g_0$. We can also contract $L$ by a
homotopy $H$ and then extend this homotopy to be $0$ on $\g_0$. Then
we find an operator $H\colon\g\colon\g[-1]$ such that
\begin{equation}\label{eqproof10}
dH+Hd=1-p
\end{equation}
Now we construct a series $\{\mathcal{F}^k\}_{k\ge 0}$ of $L_\infty$
morphisms from $\g$ to $\g$. The limit one $\mathcal{F}^{\infty}$
will be $\mathcal{F}$. We set $\mathcal{F}^0=\Id_{\g}$. Each
$\mathcal{F}^i$ is homotopical equivalent to the previous one
$\mathcal{F}^{i-1}$ in the sense of Section 2. Namely, we set
\begin{equation}\label{eqproof11}
\mathcal{F}^{k+1}=\mathcal{F}^k_{X_k}
\end{equation}
where
\begin{equation}\label{eqproof12}
X_k=-H\circ \mathcal{F}^k
\end{equation}
More precisely,
\begin{equation}\label{eqproof13}
\mathcal{F}^{k+1}=\exp(\ad (-H\circ
\mathcal{F}^k))(\mathcal{F}^k)+f(\ad (-H\circ
\mathcal{F}^k))(d_\Bbbk \mathcal{F}^k)
\end{equation}
where $f(z)=(\exp(z)-1)/z$.

\begin{lemma}
The Taylor components $(\mathcal{F}^k)_i$ of the $L_\infty$ morphism
$\mathcal{F}^k\colon\g\to\g$ take values in $\g_0=p(\g)$ for $i\le
k$. The Taylor components $(\mathcal{F}^{k+1})_i$ and
$(\mathcal{F}^k)_i$ coincide for $i\le k$. That is, the sequence of
$L_\infty$ morphisms $\{\mathcal{F}^k\}$ stabilizes, and the limit
$L_\infty$ morphism $\mathcal{F}^\infty\colon\g\to\g$ takes values
in $\g_0$.
\begin{proof}
We compute directly from (\ref{eq2.6}):
$(\mathcal{F}^1)_1=\Id-(\Id-p)=p$, that is, the claim is true for
$k=1$.

Suppose that all $i$-th Taylor components of $\F^k$ take values in
$\g_0$ for all $i\le k$. Then by definition, the $L_\infty$ morphism
$\F^{k+1}$ is obtained from $\F^k$ by the action of the homotopy
$-H\circ\F^k$. Then all $i$-th Taylor components of $-H\circ\F^k$
are 0 for $i\le k$, because $H\circ p=0$ (we constructed the
homotopy to be 0 on $\g_0$). This proves that $\F^k$ and $\F^{k+1}$
coincide up to the $k$-th Taylor component. It remains to prove that
the $(k+1)$-st Taylor component of $\F^{k+1}$ takes values in
$\g_0$. Indeed, one has:
\begin{equation}\label{tamarkin6}
(\F^{k+1})_{k+1}(x_1,\dots,x_{k+1})=(\F_k)_{k+1}(x_1,\dots,x_{k+1})-[d,(H\F^k)_{k+1}(x_1,\dots,x_{k+1})]
\end{equation}
The second summand is
\begin{equation}\label{tamarkin7}
[d,(H\F^k)_{k+1}(x_1,\dots,x_k)]=[d,H](\F^k)_{k+1}(x_1,\dots,x_{k+1}))-H[d,(\F^k)_{k+1}](x_1,\dots,x_{k+1})
\end{equation}
But
\begin{equation}\label{tamarkin8}
\begin{aligned}
\ &[d,(\F^k)_{k+1}](x_1,\dots,x_{k+1})=\\ &\pm
\sum_{i<j}\pm(\F^k)_k([x_i,x_j],x_1,\dots,\hat{x}_i,\dots,\hat{x}_j,\dots,x_{k+1})\\
+&\frac12\sum_{a+b=k+1}\frac1{a!}\frac1{b!}\pm[(\F^k)_a(x_{\sigma_1},\dots),(\F^k)_b(x_{\sigma_{a+1}},\dots)]
\end{aligned}
\end{equation}
takes values in $\g_0$ by induction, and because $\g_0$ is a Lie
subalgebra (is closed under the bracket). Therefore, the second
summand in (\ref{tamarkin7}) is 0 because $H|_{\g_0}=0$. We have:
\begin{equation}\label{tamarkin9}
(\F^{k+1})_{k+1}=(\F^k)_{k+1}-[d,H](\F^k)_{k+1}=p(\F^k)_{k+1}
\end{equation}

Lemma is proven.

\end{proof}
\end{lemma}
Now we prove that the compositions $\mathcal{F}^\infty\circ i\sim
\Id_{\g_0}$ and $i\circ \mathcal{F}^\infty\sim\Id_\g$, where $\sim$
stands for the homotopical equivalence of $L_\infty$ morphisms.

The second relation is trivial, because $i\circ
\mathcal{F}^\infty=\mathcal{F}^\infty$ which is homotopical
equivalent to $\mathcal{F}^0=\Id$ by the construction. (We use the
Campbell-Baker-Hausdorff formula).

Let us prove the first one. When all $x_i$ belong to $\g_0$, one
easily sees that
$(\mathcal{F}^1)_k(x_1,\dots,x_k)=(\mathcal{F}^0)_k(x_1,\dots,x_k)$
for any $k\ge 1$. It proves that
$(\mathcal{F}^\infty)_k(x_1,\dots,x_k)=(\mathcal{F}^0)_k(x_1,\dots,
x_k)$, therefore, $\mathcal{F}^\infty\circ i=\Id_{\g_0}$.

Proposition is proven.
\end{proof}
\end{proposition}

\subsection{We end to prove Theorem 1}
We are now return to the notations of Section 4.1. By Proposition
4.4, the maps $\mathcal{L}(\mathcal{M}_i)$ (where
$\mathcal{M}_i\colon H^\mb(\g_i)\to \g_i$, $(i=1,2)$, are the maps
obtained by the planar graph construction), being quasi-isomorphic
imbedding of dg Lie algebras, are homotopically invertible. Then by
Lemma 4.3 the maps $\mathcal{M}_i$ are homotopically invertible.
(Lemma 4.3 is applicable because $H^\mb(\g_i)$ is an $L_\infty$
algebra and $\g_i$ is a dg Lie algebra). Consider now the diagram
(\ref{eq4.1}). Then the dotted arrow exists up to homotopy, and
induces an {\it isomorphism} in the first Taylor component, by Lemma
4.1. Then it can be inverted by the inverse function theorem. We
conclude that the quasi-isomorphism $\mathcal{F}\colon\g_1\to\g_2$
can be homotopically inverted.

\qed

\section{The category $\mathcal{H}om_{dg}$ is a localization} We have proved in the
previous Section that in the category $\mathcal{H}om_{dg}$ the
invertible morphisms are the quasi-isomorphisms. That is, the first
property of the localization of category from Section 1.1 is
satisfied. Here we prove the second property, that is, that the
category $\mathcal{H}om_{dg}$ is a localization of the category
$\mathcal{L}ie_{dg}$ by quasi-isomorphisms (and then it is a
localization of the category of dg Lie algebras and dg Lie maps by
quasi-isomorphisms). This is the following result:
\begin{theorem}
Suppose $F\colon \mathcal{L}ie_{dg}\to\mathfrak{X}$ is a functor
which maps any quasi-isomorphism to an invertible morphism. Then $F$
maps any two morphisms
$\varphi_1,\varphi_2\in\mathrm{Mor}(\g_1,\g_2)$ homotopic in the
sense of Section 2 to the same morphism $F(\varphi_1)=F(\varphi_2)$
in $\mathfrak{X}$.
\end{theorem}
It follows from the definition of the localization of a category
(see Section 1.1) that any two localizations are equivalent.
Therefore, the category $\mathcal{H}om_{dg}$ is equivalent to the
corresponding Quillen categories. Also it follows from the Theorem
above that our definition of a homotopy coincides with the Quillen's
one.

We prove Theorem 2 in the rest of this Section.

\subsection{A cylinder construction: the Lie case}
Let $\g$ be a dg Lie algebra.  Define a dg Lie algebra $\Cyl(\g)$ as
follows. As a graded vector space, it is the direct sum
$\Cyl(\g)=\Omega^0([0,1],\g)\oplus \Omega^1([0,1],\g)$. Here
$\Omega^0([0,1],\g)$ is the space of smooth functions on the
interval $[0,1]$ with values in $\g$, and $\Omega^1([0,1],\g)$ is
the space of smooth 1-forms on $[0,1]$ with values in $\g$. The
grading is $\deg(f(t)\cdot g)=\deg g$ for a smooth function $f(t)$
on $[0,1]$ and a homogeneous $g\in\g$, and $\deg(f(t)dt\cdot g)=\deg
g+1$. Introduce a differential and a Lie bracket on $\Cyl(\g)$. The
differential is the de Rham differential. It maps $f(t)\cdot g$ to
$f^\prime(t)dt\cdot g$. The bracket is linear over smooth functions
on the interval and comes from the Lie bracket on $\g$. In
particular, $[\Omega^0([0,1],\g),\Omega^0([0,1],\g)]\subset
\Omega^0([0,1],\g)$, $[\Omega^0([0,1],\g),\Omega^1([0,1],\g)]\subset
\Omega^1([0,1],\g)$, and
$[\Omega^1([0,1],\g),\Omega^1([0,1],\g)]=0$.

It is clear that this makes $\Cyl(\g)$ a dg Lie algebra.

This Lie algebra is quasi-isomorphic to $\g$, this fact can be
considered as the Poincar\'{e} lemma applied to the interval
$[0,1]$.

For any $s\in [0,1]$ one has a map $p_s\colon\Cyl(\g)\to\g$ of dg
Lie algebras, defined as follows: it is 0 on $\Omega^1([0,1],\g)$
and is the value for $t=s$ (the specialization) on
$\Omega^0([0,1],\g)$.

\begin{lemma}
Let $F\colon\mathcal{L}ie_{dg}\to\mathfrak{X}$ be any functor which
maps $L_\infty$ quasi-isomorphisms to invertible morphisms, and let
$s_1,s_2\in [0,1]$ be any two points. Then $F(p_{s_1})=F(p_{s_2})$.
(It is enough if the functor $F$ maps quasi-isomorphisms of dg Lie
algebras to invertible morphisms).
\begin{proof}
Define a map $\sigma\colon \g\to\Cyl(\g)$ which maps $g\mapsto
(g,0)$ to the constant function equal to $g$ and zero 1-form. It is
clear that $\sigma$ is a map of dg Lie algebras.

Now for any $s\in [0,1]$ we have: $p_s\circ\sigma=\Id_{\g}$.
Applying the functor $F$, we get: $F(p_s)\circ
F(\sigma)=\Id_{F(\g)}$. The map $\sigma$ is a quasi-isomorphism of
dg Lie algebras, and, therefore, $F(\sigma)$ is invertible morphism.
Then for any $s\in[0,1]$ we have: $F(p_s)=(F(\sigma))^{-1}$.
\end{proof}
\end{lemma}

\subsection{}
Let $\g_1,\g_2$ be dg Lie algebras, and let $U_0\colon \g_1\to \g_2$
be an $L_\infty$ map. Then $U_0$ can be considered as a solution of
the Maurer-Cartan equation in $(\Bbbk(\g_1,\g_2))^1$. Let
$H\in(\Bbbk(\g_1,\g_2)^0$ be a "homotopy", and let an $L_\infty$
morphism $U_1\colon\g_1\to\g_2$ is the homotopic $L_\infty$
morphism, given by the formula (see Section 2):

\begin{equation}\label{eq5.2}
U_1=\exp(ad(H))(U_0)+f(H)(d_\Bbbk H)
\end{equation}
where $f(z)=(\exp(z)-1)/z$.

Then we have the following lemma:
\begin{lemma}
There exists an $L_\infty$ map $U_\Cyl\colon \g_1\to\Cyl(\g_2)$ such
that the diagrams
\begin{equation}\label{eq5.1}
\xymatrix{ &&\Cyl(\g_2)\ar[dd]^{p_i}\\
\\
\g_1\ar[rruu]^{U_\Cyl}\ar[rr]^{U_i}&&\g_2}
\end{equation}
($i=0,1$) are commutative, where $p_0$ and $p_1$ are the maps
$p_s\colon\Cyl(\g_2)\to\g_2$ for $s=0$ and $s=1$.
\end{lemma}
First of all, we have the following corollary:
\begin{corollary}
In the notations of Lemma 5.2 one has $F(U_0)=F(U_1)$.
\begin{proof}
From the commutative diagrams (\ref{eq5.1}) one has:
$F(U_0)=F(p_0)\circ F(U_\Cyl)$ and $F(U_1)=F(p_1)\circ F(U_\Cyl)$.
But $F(p_0)=F(p_1)$ by Lemma 5.2.
\end{proof}
\end{corollary}
This Corollary is Theorem 3 for $\mathcal{H}om_{dg}$.
\begin{remark}
The lemma gives an explicit construction of the "right cylinder" in
the category $\mathcal{L}ie_{dg}$ of dg Lie algebras and $L_\infty$
maps between them. In the category of topological spaces, the right
cylinder is constructed as follows. Consider the space $Y^I$ of all
free paths $\alpha\colon I\to Y$ where $I=[0,1]$ is the unit
interval. This space is the topological right cylinder of $Y$. There
are two projections $p_0,p_1\colon Y^I\to Y$ which are the values in
the end-points of the interval. Let $f_0,f_1\colon X\to Y$ are two
homotopic maps between topological spaces (the homotopy is
understood in the usual, left, sense). Then clearly there exists a
map $h\colon X\to Y^I$ such that the diagrams
$$
\xymatrix{ & Y^I\ar[d]^{p_i}\\
X\ar[ru]^{h}\ar[r]^{f_i}&Y}
$$
$(i=0,1)$ are commutative. This $h$ can be constructed directly from
the left homotopy between $f_0$ and $f_1$. Thus, for the category of
topological spaces the left and the right homotopies are equivalent.
\end{remark}

\subsection{A proof of Lemma 5.2}
We construct an $L_\infty$ map $U_\Cyl\colon\g_1\to\Cyl(\g_2)$ with
the demanded properties.

Recall from Section 2, that for any $H\in\Bbbk(\g_1,\g_2)^1$ one has
a vector field
\begin{equation}\label{eq5.3.1}
\frac{dF}{dt}=-dH+[H,F]
\end{equation}
on $\Bbbk(\g_1,\g_2)^1$ tangent to the Maurer-Cartan quadric. It can
be integrated to a solution $F\colon [0,1]\to \Bbbk(\g_1,\g_2)^1$
when we fix the initial condition $F(0)=U_0$. When the initial
condition belongs to the Maurer-Cartan quadric, $F(t)$ belongs to
this quadric for all $t\in[0,1]$.

The solution $F(t)$ is given by the formula
\begin{equation}\label{eq5.3.2}
F(t)=\exp(ad(tH))(F(0))+f(tH)(d_\Bbbk H)
\end{equation}
where $f(z)=(\exp(z)-1)/z$. It follows from the fact that
(\ref{eq5.3.1}) can be rewritten as
\begin{equation}\label{eq5.3.3}
\frac{d(F+d_\Bbbk)}{dt}=[H,F+d_\Bbbk]
\end{equation}
We define the $L_\infty$ map $U_\Cyl\colon\g_1\to\Cyl(\g_2)$ as an
element in $\Bbbk(\g_1,\Cyl(\g_2))^1$ defined by
\begin{equation}
\theta\mapsto (F(t), H\cdot dt)
\end{equation}
where $\theta\in C_+(\g_1,\mathbb{C})$. Here
$F(t)\in\Omega^0([0,1],\g_2)$ and $H\cdot dt\in
\Omega^1([0,1],\g_2)$. We need to prove that this element
$U_\Cyl\in\Bbbk(\g_1,\Cyl(\g_2))^1$ satisfies the Maurer-Cartan
equation (that is, defines an $L_\infty$ morphism from $\g_1$ to
$\Cyl(\g_2)$.

For this denote $\alpha=F(t)+H\cdot dt\in\Bbbk(\g_1,\Cyl(\g_2))^1$.
The Maurer-Cartan equation is
\begin{equation}\label{eq5.3.4}
d_\Bbbk\alpha +\frac12[\alpha,\alpha]=0
\end{equation}
This equation is equivalent to two equations, correspondind to the
two components in $\Cyl(\g_2)=\Omega^0\oplus\Omega^1$.

For the $\Omega^0$ component the Maurer-Cartan equation
(\ref{eq5.3.4}) is equivalent to the fact that for any $t\in[0,1]$
the element $F(t)\in\Bbbk(\g_1,\g_2)$ satisfies the Maurer-Cartan
equation, which is clear by the construction.

The $\Omega^1$ component of (\ref{eq5.3.4}) is
\begin{equation}\label{eq5.3.5}
d_{DR}F(t)+d_\Bbbk(H)dt+[Hdt,F(t)]=0
\end{equation}
which is equivalent to (\ref{eq5.3.1}). (Note that $[H\cdot
dt,F(t)]=-[H,F]\cdot dt$ by (\ref{signs}), because $\deg F(t)=1$ and
$\deg dt=1$).

The $L_\infty$ map $U_\Cyl\colon\g_1\to\Cyl(\g_2)$ is constructed.
One should check that the diagrams (\ref{eq5.1}) are commutative.
The commutativity for $p_0$ follows from the construction, and the
commutativity for $p_1$ follows from $F(1)=U_1$ if $F(0)=U_0$.

Lemma 5.2 is proven. \qed

Theorem 2 is proven. \qed

\bigskip
\bigskip

Faculty of Science, Technology and Communication, Campus
Limpertsberg, University of Luxembourg,
162A avenue de la Faiencerie, L-1511 LUXEMBOURG\\
{\it e-mail}: {\tt borya$\_$port@yahoo.com}

\end{document}